\documentclass[12pt]{article}
\usepackage{amsmath, amssymb, amsthm, amsfonts}
\usepackage{indentfirst}
\usepackage[mathscr]{euscript}
\usepackage{amscd}
\usepackage{color}

\numberwithin{equation}{section}

\theoremstyle{plain}
\newtheorem{theorem}{Theorem}\numberwithin{theorem}{section}
\newtheorem{lemma}{Lemma}\numberwithin{lemma}{section}
\newtheorem{proposition}{Proposition}\numberwithin{proposition}{section}
\newtheorem{corollary}{Corollary}\numberwithin{corollary}{section}
\theoremstyle{definition}
\newtheorem{definition}{Definition}\numberwithin{definition}{section}
\theoremstyle{remark}
\newtheorem{remark}{Remark}\numberwithin{remark}{section}

\newcommand{\R}{\mathbb{R}}

\newcommand{\s}{\mathbb{S}}

\newcommand{\F}{\mathcal{F}}
\newcommand{\bx}{\mathbf{x}}
\newcommand{\by}{\mathbf{y}}

\title{Compactly supported reproducing kernels\\ for $L^2$-based Sobolev spaces\\
and Hankel-Schoenberg transforms}

\author{Yong-Kum Cho\thanks{Department of Mathematics, College of Natural Sciences, Chung-Ang University,
84 Heukseok-Ro, Dongjak-Gu, Seoul 156-756, Korea (e-mail: ykcho@cau.ac.kr)}
\footnote{This research was supported by National Research Foundation of Korea Grant
funded by the Korean Government (\# 20160925).}}

\begin{document}

\maketitle

\bigskip

\begin{itemize}
\item[{}] {\bf Abstract.} We exhibit three classes of compactly supported
functions which provide reproducing kernels for the Sobolev spaces
$H^\delta(\R^d)$ of arbitrary order $\,\delta>d/2.\,$
Our method of construction is based on a new class of oscillatory integral transforms that
incorporate radial Fourier transforms and Hankel transforms.
\end{itemize}

\bigskip
{\small
\begin{itemize}
\item[{}]{\bf Keywords.} Askey's class, Bessel function, Bessel potential kernel, binomial density,
Fourier transform, generalized hypergeometric function, Hankel-Schoenberg transform, Hilbert space, positive definite,
reproducing kernel, Sobolev space, Wendland's function.

\item[{}] 2010 Mathematics Subject Classification: 33C10, 41A05, 42B10, 60E10.
\end{itemize}}

\newpage

\section{Introduction}
In this paper we shall deal with the problem of constructing {\it compactly supported radial} functions $\Phi$
on $\R^d$ such that the symmetric kernels $\Phi(\bx-\by)$ could serve as reproducing kernels for Sobolev spaces
under appropriate inner products. Due to advantageous aspects in practical applications,
the problem has become an important issue in various fields of Mathematics including the theory of
interpolations, spatial statistics and machine learning.

In their pioneering work \cite{AK}, N. Aronszajn and K. T. Smith introduced the Sobolev space
$\,H^\delta(\R^d)\,$ of order $\,\delta>0\,$ as the space of Bessel potentials defined by
convolutions $\,(G_{\delta/2}\ast u)(\bx),\,$ where $\,u\in L^2(\R^d)\,$ and
\begin{equation}\label{G1}
G_\delta(\bx) = \frac{1}{2^{\delta-1 + \frac d2}\,\pi^{\frac d2}\,
\Gamma(\delta)}\,K_{\delta-\frac d2}(|\bx|) |\bx|^{\delta - \frac d2}.
\end{equation}
As usual, $\,|\bx| = \sqrt{\bx\cdot\bx}\,$ denotes the Euclidean norm for each $\,\bx\in\R^d\,$
and $K_{\delta- d/2}$ stands for the modified Bessel function of order
$\,\delta-d/2.$

Often referred to as Mat\'ern functions (see e.g. \cite{G1}), the Bessel potential kernels $G_\delta$ are
are integrable with the Fourier transforms
\begin{align}\label{G2}
\widehat{G_\delta}(\xi) =\int_{\R^d} e^{-i \xi\cdot\bx}\,G_\delta(\bx) d\bx= \left( 1+|\xi|^2\right)^{-\delta}.
\end{align}
As a consequence, the Sobolev space of order $\delta$ may be identified with
\begin{equation*}
H^\delta(\R^d) =\left\{ u\in L^2(\R^d) : \left( 1+|\cdot|^2\right)^{\delta/2}\,\widehat{u}\in L^2(\R^d)\right\},
\end{equation*}
which becomes a Hilbert space under the inner product
\begin{align*}
\bigl(u,\,v\bigr)_{H^\delta(\R^d)} = (2\pi)^{-d}\int_{\R^d}  \left( 1+|\xi|^2\right)^{\delta}\widehat{u}(\mathbf{\xi})
\overline{\,\widehat{v}(\mathbf{\xi})} \,d\mathbf{\xi}.
\end{align*}

In the case $\,\delta> d/2,\,$  N. Aronszajn and K. T. Smith noticed further that
$\,H^\delta(\R^d)\subset C(\R^d)\,$ continuously and $H^\delta(\R^d)$ is a
reproducing kernel Hilbert space with kernel $G_{\delta}(\bx - \by)$, that is,
for every $\,u\in H^\delta(\R^d),\,\bx\in\R^d,\,$
\begin{align*}
&{\rm(i)}\quad G_{\delta}(\cdot-\bx)\in H^\delta(\R^d)\quad\text{and}\\
&{\rm(ii)}\quad u(\bx) = \big(u, \,G_\delta(\cdot-\bx)\big)_{H^\delta(\R^d)}
\end{align*}
(we also refer to A. P. Caler\'on \cite{C} and the appendix for a brief additional description of the
Bessel potential kernels $G_\delta$).

In connection with the problem of our consideration, there is a standard framework on reproducing kernel Hilbert spaces
of functions on $\R^d$ which resembles the structure of Sobolev spaces and reads as follows.
For a given real-valued {\it positive definite}
function $\,\Phi\in C(\R^d)\cap L^1(\R^d),\,$
if we define
\begin{align*}
& \mathcal{F}_\Phi(\R^d) = \left\{ u\in C(\R^d)\cap L^2(\R^d) :
\int_{\R^d} \big|\widehat u(\mathbf{\xi})\big|^2\frac{d\xi}{\,\widehat{\Phi}(\mathbf{\xi})\,}<\infty\right\},\\
&\qquad\bigl(u,\,v\bigr)_{\mathcal{F}_\Phi(\R^d)} = (2\pi)^{-d}\int_{\R^d} \widehat{u}(\mathbf{\xi})\,
\overline{\widehat{v}(\xi)}\,\frac{d\xi}{\,\widehat{\Phi}(\mathbf{\xi})\,},
\end{align*}
then $\mathcal{F}_\Phi(\R^d)$ becomes a Hilbert space with a reproducing kernel $\Phi(\mathbf{x}-\mathbf{y})$
(see \cite{S1}, \cite{We} and also \cite{Ar} for more general properties).

On account of this framework, we shall focus on constructing
compactly supported radial functions $\,\Phi\in C(\R^d) \cap L^1(\R^d)\,$
which are positive definite and subject to the Fourier transform estimates
\begin{equation}\label{G5}
C_1 (1+|\xi|^2)^{-\delta}\le  \widehat\Phi(\xi)\le C_2 (1+|\xi|^2)^{-\delta}
\end{equation}
for some $\,\delta> d/2\,$ and for some positive constants $\,C_1, C_2.$

An initiative construction had been started by H. Wendland (\cite{We1}, \cite{We}) who introduced
a family of polynomials on $[0, \infty)$ defined by
\begin{align}\label{G3}
\left\{\begin{aligned} {P_{d, m}(r)} &{\,\,= c_m\int_{r}^1 \left(t^2-r^2\right)^{m-1} t
(1-t)^{\left[\frac d2\right] + m +1} dt,}\\
{P_{d}(r)} &{\,\,= (1- r)^{\left[\frac d2\right] +1}\,,}\end{aligned}\right.
\end{align}
for $\,0\le r\le 1\,$ and zero otherwise,
where $m$ is a positive integer and $c_m$ is a constant,
and proved $P_{d, m}(|\bx-\by|)$ is a reproducing kernel for the Sobolev space
$\,H^{\frac{d+1}{2} + m}(\R^d)\,$ and so is $P_{d}(|\bx-\by|)$ for
$\,H^{\frac{d+1}{2}}(\R^d)\,$ if $\,d\ge 3.$

In an attempt to cover the missing cases,
R. Schaback (\cite{S2}) introduced a family of non-polynomial functions defined by
\begin{equation}\label{G4}
R_{d, m}(r) = c_m\int_{r}^1 \left(t^2-r^2\right)^{m- \frac 12} t
(1-t)^{\frac d2 + m +1} dt
\end{equation}
for $\,0\le r\le 1\,$ and zero otherwise,
where $m$ is a nonnegative integer, and proved
$R_{d, m}(|\bx-\by|)$ is a reproducing kernel for
$\,H^{\frac d2 + m + 1}(\R^d)\,$ if $d$ is even (see S. Hubbert \cite{H} for computational aspects).

In order to deal with fractional orders, A. Chernih and S. Hubbert (\cite{CH}) further generalized Wendland's functions
in the form
\begin{equation}\label{G5}
S_{d, \alpha}(r) = c_\alpha \int_{r}^1 \left(t^2-r^2\right)^{\alpha-1} t(1-t)^{\frac {d+1}{2} +\alpha} dt
\end{equation}
for $\,0\le r\le 1\,$ and zero otherwise, where $\,\alpha>0\,$ and $c_\alpha$ is a constant, and proved that
$S_{d, \alpha}(|\bx-\by|)$ is a reproducing kernel for $\,H^{\frac {d+1}{2} +\alpha}(\R^d).\,$

Our primary aim in the present paper is to obtain a family of compactly supported
radial functions which provide reproducing kernels for
the Sobolev spaces $H^\delta(\R^d)$ of any order $\,\delta>d/2\,$ in a unified manner
and thereby cover all of the missing cases left open in this subject.

The method of our construction will be based on a new class of oscillatory integral transforms, to be called
Hankel-Schoenberg transforms hereafter, which incorporate Fourier transforms of radial functions and
classical Hankel transforms. Apparently useful in any situation where Fourier transforms of radial functions are involved,
our secondary purpose is to bring Hankel-Schoenberg transforms to attention and establish their basic properties.

In consideration of Euler's binomial densities as possible candidates, we shall begin with evaluating
their Hankel-Schoenberg transforms in terms of generalized hypergeometric functions whose asymptotic behaviors
are well investigated, e.g., by Y. L. Luke \cite{L}. We then select those
binomial densities whose Hankel-Schoenberg transforms are strictly positive by the criteria of J. Fields and
M. Ismail \cite{FI} and apply a continuous version of dimension walks to obtain the desired classes of functions.

As it will be presented in detail, we shall exhibit three different classes of compactly supported
functions which provide reproducing kernels for the Sobolev spaces $H^\delta(\R^d)$
of order
$$\,\delta = \frac{d+1}{2}, \quad\delta>\max\,\left(1, \,\frac d2\right), \quad \delta> \frac d2\,,$$
separately. One of these classes include the compactly supported functions of Wendland, Schaback,
Chernih and Hubbert as special instances.

A distinctive feature of our construction is that the Fourier transform
is explicit, which enables us to specify the inner product about which the reproducing property holds.
As an illustration, it will be shown that the function
$\,A_2(x) = (1 - |x|)_+^2\,, x\in\R,\,$ has the Fourier transform
$$\widehat{A_2}(\xi) = \frac{4}{\xi^2}\left( 1-
\frac{\sin \xi}{\xi}\right),$$
which is strictly positive and behaves like the Cauchy-Poisson kernel,
and $A_2(x-y)$ is a reproducing kernel for $H^1(\R)$ under the inner product
$$\bigl(u,\,v\bigr)_{A_2(\R)} = \frac{1}{8\pi}\,\int_{-\infty}^\infty \widehat{u}(\mathbf{\xi})\,
\overline{\widehat{v}(\xi)}\,\frac{\xi^3\,d\xi}{\,\xi - \sin \xi\,}.$$

\medskip

\paragraph{Notation.} We shall use the following notation in what follows.
\begin{itemize}
\item The Euler beta function will be denoted by
$$B(a, \,b) = \int_0^1 t^{a-1} (1-t)^{b-1} dt\qquad(a>0, \,b>0).$$
\item The generalized hypergeometric functions
will be denoted by
\begin{equation*}
{}_pF_q\left(a_1, \cdots, a_p;\,b_1, \cdots, b_q;\,r\right)
=\sum_{k=0}^\infty\frac{\left(a_1\right)_k\cdots\left(a_p\right)_k}{k!
\left(b_1\right)_k\cdots\left(b_q\right)_k}\,r^k
\end{equation*}
in which $\,(a)_k = a(a+1)\cdots (a+k-1)\,$ if $\,k\ge 1\,$ and $\,(a)_0 = 1\,$
for any real number $a$.
\item The positive part of $\,x\in\R\,$ will be denoted by $\,x_+ = \max (x,\,0).$
\item We shall write $\,f(x) \approx g(x)\,$ for $\,x\in X\,$ for two real-valued functions
$\,f, g\,$ defined on $X$ to indicate there exist positive constants $\,c_1, c_2\,$ such that
$\,c_1\, g(x) \le f(x) \le c_2\, g(x)\,$ for all $\,x\in X.$
\end{itemize}

\bigskip

\section{Positive definite functions}
We recall that a function $\Phi$ on $\R^d$ is said to be {\it positive semi-definite}
if
\begin{align*}
\sum_{j=1}^N\sum_{k=1}^N \,\Phi\left(\bx_j - \bx_k\right) z_j \overline{z_k}\ge 0
\end{align*}
for any choice of $\,z_1, \cdots, z_N\in\mathbb{C}\,$ and $\,
\bx_1, \cdots, \bx_N\in\R^d.\,$ If equality holds only when
$\,z_1=\cdots = z_N=0,\,$ $\Phi$ is said to be {\it positive definite}.

A well-known theorem of S. Bochner states that a continuous function $\Phi$ is
positive semi-definite if and only if it is the Fourier transform of some finite
nonnegative Borel measure $\mu$ on $\R^d$. If the carrier of $\mu$
contains an open set, then $\,\Phi = \widehat{\mu}\,$ is positive definite.
In particular, the Fourier transform of a nonnegative function $\,f\in L^1(\R^d)\,$
is positive definite if the essential support of $f$ contains an open set (see \cite{We})
\footnote{The carrier of a nonnegative Borel measure $\mu$ on $\R^d$ is defined to be
$$\R^d\setminus \left\{O\subset\R^d : O \,\,\text{is open and}\,\,
\mu(O) = 0\right\}.$$
If $\mu$ is absolutely continuous with respect to Lebesgue measure, $\,d\mu(\bx) = f(\bx)\, d\bx\,$ with
a nonnegative $\,f\in L^1(\R^d)\,,$ then the carrier of $\mu$ equals to the essential support of $f$,
the complement of the largest open subset of $\R^d$ on which $\,f = 0\,$ almost everywhere.}

A univariate function $\phi$ on $[0, \infty)$ is said to be
{\it positive semi-definite or positive definite on $\R^d$} if the radial extension
$\,\bx\mapsto \phi(|\bx|), \,\bx\in\R^d,\,$ is positive semi-definite or positive definite
in the above sense. To state sufficient or necessary conditions in terms of Fourier transforms,
we shall introduce the following kernels, more extensive than being needed, which will
serve as the kernels of Hankel-Schoenberg transforms to be studied later.

\smallskip

\begin{definition}\label{def1}
For $\,\lambda>-1,\,$ define $\,\Omega_\lambda : \R\to \R\,$ by
\begin{align*}
\Omega_\lambda(t) &=
\Gamma(\lambda+1)\sum_{k=0}^\infty\frac{(-1)^k}{k!\,\Gamma(\lambda +k +1)}\,\left(\frac t 2\right)^{2k}\\
&=\Gamma(\lambda+1)\left(t/2\right)^{-\lambda} J_\lambda(t),
\end{align*}
where $J_\lambda$ denotes the Bessel function of the first kind of order $\lambda$.
\end{definition}

\smallskip
In the special case $\,\lambda = (d-2)/2\,,$ with $\,d\ge 2\,$ a positive integer, $\Omega_\lambda$
arises on consideration of the Fourier transform of the area measure $\sigma$ on the unit sphere
$\s^{d-1}$ of the Euclidean space $\R^d$ in the form
\begin{equation*}
\frac{1}{\left|\s^{d-1}\right|}\, \int_{\s^{d-1}} e^{-i \xi\cdot \bx}\,d\sigma(\bx)
=\Omega_{\frac{d-2}{2}}(|\xi|).
\end{equation*}

An immediate consequence is that if
$F$ is integrable and radial with $\,F(\bx) = f(|\bx|)\,$ for some univariate function $f$
on $[0, \infty)$, then its Fourier transform is easily evaluated as
\begin{align*}
\widehat{F}(\xi) &= \int_0^\infty \left(\int_{\s^{d-1}} e^{-i t\xi\cdot \bx}\,d\sigma(\bx)
\right) f(t) t^{d-1} dt\\
&= |\s^{d-1}|\int_0^\infty \Omega_{\frac{d-2}{2}}(|\xi|t) f(t) t^{d-1} dt.
\end{align*}
Since it is simple to find $\,\Omega_{-1/2}(t) =\cos t\,$ by definition, this formula continues
to hold true for $\,d=1\,$ if we interpret $\,|\s^{0}| =2.$

In summary, we have the following which are substantially due to I. J. Schoenberg \cite{Sc}
(see also  \cite{GMO}, \cite{MS}, \cite{Stw}, \cite{We}).

\smallskip

\begin{proposition}\label{fourier} For $\,f\in L^1\left([0, \infty), \,t^{d-1} dt\right),\,$
put
\begin{equation*}
\phi(r) = \int_0^\infty\Omega_{\frac{d-2}{2}}\left(rt\right) f(t) t^{d-1} dt\qquad(r\ge 0).
\end{equation*}
\begin{itemize}
\item[\rm(i)] If $\,F(\bx) = f(|\mathbf{x}|),\,\mathbf{x}\in\R^d,\,$ then the Fourier transform of $F$ is given by
\begin{equation*}\label{F2}
\widehat{F} (\xi) = \frac{2\pi^{d/2}}{\Gamma(d/2)}\,\phi(|\mathbf{\xi}|).
\end{equation*}
\item[\rm(ii)] If $f$ is nonnegative and the essential support of $f$ contains an open interval,
then $\phi$ is positive definite on $\R^d$.
\end{itemize}
\end{proposition}

\smallskip

\begin{proof}
Part (i) is what we have mentioned as above.
Concerning part (ii), if the essential support of a nonnegative function
$f$ contains an open interval, say, $\,I= (a, b),\,$
then the essential support of $\,F(\bx) = f(|\bx|)\,$ contains the open annulus
$\,\{\bx\in\R^d: a<|\bx|<b\}\,$ and the assertion follows.
\end{proof}

\smallskip
\begin{remark}
The integral defined in the statement is often called the $d$-dimensional radial Fourier transform
and formally denoted as
\begin{equation}\label{P2}
\F_d(f)(r) = \frac{2\pi^{d/2}}{\Gamma(d/2)}\,
\int_0^\infty\Omega_{\frac{d-2}{2}}\left(rt\right) f(t) t^{d-1} dt\qquad(r\ge 0).
\end{equation}
As the kernel $\Omega_{\frac{d-2}{2}}$ will be shown to be uniformly bounded,
the integral makes sense on the class of finite Borel measures on $[0, \infty)$.
Indeed, Schoenberg's original theorem states that a continuous function $\phi$ on $[0, \infty)$
is positive semi-definite on $\R^d$ if and only if
$$\phi(r) = \int_0^\infty\Omega_{\frac{d-2}{2}}\left(rt\right) d\nu(t)\qquad(r\ge 0)$$
for some finite nonnegative Borel measure $\nu$ on $[0, \infty)$.
\end{remark}

\bigskip

\section{Hankel-Schoenberg transforms}
As it is classical (see \cite{Wa} for instance), the
Hankel transforms of a function $\,f\in L^1\left([0, \infty),\,\sqrt t dt\right)\,$ refer to the integrals of type
$$\int_0^\infty J_\lambda(rt) f(t) t dt\qquad(\lambda\ge -1/2).$$

As a generalization of both Fourier transforms of radial functions and Hankel transforms,
we shall consider the following integral transforms.

\smallskip

\begin{definition}
The Hankel-Schoenberg transform of order $\,\lambda\ge -1/2\,$
of a Lebesgue measurable function $f$ on $[0, \infty)$ is defined to be
\begin{equation*}
\phi(r) = \int_0^\infty\Omega_\lambda(rt) f(t) dt \qquad(r\ge 0)
\end{equation*}
whenever the integral on the right side converges.
\end{definition}

\smallskip

The definition indeed makes sense under various conditions on $f$. For this matter,
we shall begin with investigating the kernels.

\medskip

\subsection{Kernels $\Omega_\lambda$}
In many aspects, each $\Omega_\lambda$ is similar in nature to the cardinal sine function
$$\frac{\sin t}{t} = \prod_{k=1}^\infty \left( 1- \frac{t^2}{k^2\pi^2}\right)\,$$
which coincides with the special case $\,\lambda = 1/2.$ To be more specific, we list
the following properties of $\Omega_\lambda$'s which are deducible from
the theory of Bessel functions $J_\lambda$ in a straightforward manner (see \cite{AS}, \cite{E}, \cite{Wa}).

\begin{itemize}
\item[(P1)] Each $\Omega_\lambda$ is of class $C^\infty(\R)$, even and uniformly bounded by $\,1=\Omega_\lambda(0).\,$
The kernels $\Omega_\lambda$ satisfy the Bessel-type differential equations
\begin{equation*}
{\Omega_\lambda}''(t) + \frac{2\lambda +1}{t}\,{\Omega_\lambda}'(t) +  \Omega_\lambda(t) = 0
\end{equation*}
as well as the Lommel-type recurrence relations
\begin{align}\label{P1}
&\qquad {\Omega_\lambda}'(t) = - \frac{t}{\,2(\lambda +1)\,}\,\Omega_{\lambda+1}(t),\nonumber\\
&\Omega_\lambda(t) -\Omega_{\lambda-1}(t) = \frac{t^2}{\,4\lambda(\lambda+1)\,}\,\Omega_{\lambda+2}(t).
\end{align}

\item[(P2)] An asymptotic formula due to Hankel states that as $\,t\to\infty,$
\begin{equation*}
\Omega_\lambda(t) = \frac{\Gamma(\lambda+1)}{\sqrt{\pi}} \left(\frac t2\right)^{-\lambda -1/2}
\left[\cos\left(t - \frac{(2\lambda+1)\pi}{4}\right) + O\left( t^{-1}\right)\right].
\end{equation*}

\item[(P3)]$\Omega_\lambda$ is oscillatory with an infinity of simple zeros.  Arranging the positive zeros of $J_\lambda$
in the ascending order $\,0<j_{\lambda, 1} < j_{\lambda, 2} < \cdots\,,$
$\Omega_\lambda$ can be represented as the infinite product
\begin{equation*}
\Omega_\lambda(t) = \prod_{k=1}^\infty
\left( 1- \frac{t^2}{j_{\lambda, k}^2}\right).
\end{equation*}

\item[(P4)] Due to Liouville, $\Omega_\lambda$ is expressible in finite terms
by algebraic and trigonometric functions if and only if $2\lambda$ is an odd integer. Indeed,
\begin{equation}
\Omega_{-1/2}(t) = \cos t\,,\quad \Omega_{1/2}(t) = \frac{\sin t}{t}
\end{equation}
and recurrence formula \eqref{P1} may be used to express $\Omega_{n + 1/2}$, with $n$ an integer,
in finite terms by elementary functions.
For example,
\begin{align}
\Omega_{3/2}(t) &= 3\left(\frac{\sin t - t\cos t}{t^3}\right),\nonumber\\
\Omega_{5/2}(t) &= 15\left[\frac{( 3- t^2)\sin t - 3t\cos t}{t^5}\right].
\end{align}

\item[(P5)] For $\,\lambda>-1/2,\,$ Poisson's integral reads
\begin{align*}
\Omega_\lambda(t) = \frac{2}{B\left(\lambda + 1/2\,,\,1/2\right)}\,
\int_0^1 \cos(t s)\, (1-s^2)^{\lambda -\frac 12}\,ds.
\end{align*}
\end{itemize}

Owing to the boundedness and asymptotic behavior of $\Omega_\lambda$ described as in (P1),
(P2), it is evident that the Hankel-Schoenberg transform of order $\lambda$ is well defined
on the class $\,L^1([0, \infty))\,$ or $\,L^1\left( [0, \infty),\,t^{-\lambda-1/2} dt\right).$

\medskip

\subsection{Inversion formula}
The Hankel-Watson inversion theorem (\cite{Wa}) states that if $\,\lambda\ge -1/2\,$
and $\,f(s)\sqrt s\,$ is integrable on $[0, \infty)$, then
\begin{equation*}
\int_0^\infty J_\lambda(rt)\left[\int_0^\infty J_\lambda(rs) f(s)s ds\right] rdr
= \frac{\,f(t+0) + f(t-0)\,}{2}
\end{equation*}
at every $\,t>0\,$ such that $f$ is of bounded variation in a neighborhood of $t$.

As it is straightforward to express Hankel-Schoenberg transforms in terms of Hankel transforms,
an obvious modification yields the following.

\smallskip

\begin{proposition}\label{inversion} {\rm (Inversion)}
For $\,\lambda\ge -1/2,\,$ assume that
\begin{equation}\label{invc1}
\int_0^\infty |f(t)| t^{-\lambda-1/2}\,dt <\infty.
\end{equation}

Then the following holds for every $\,t>0\,$ at which $f$ is continuous:
\begin{align*}
\left\{\aligned &{\phi(r) = \int_0^\infty \Omega_\lambda(rt) f(t) dt\quad(r>0)\quad \text{implies}}\\
 &{f(t) = \frac{t^{2\lambda+1}}{4^\lambda\left[\Gamma(\lambda+1)\right]^2}\,
 \int_0^\infty \Omega_\lambda(rt)\, \phi(r) r^{2\lambda+1}dr.}\endaligned\right.
\end{align*}
\end{proposition}

\smallskip

\begin{remark}
In the case $\,\lambda = (d-2)/2\,,$ this formula may be considered
as an alternative of the Fourier inversion theorem for radial functions.
Useful to the present circumstance is the inversion of $\,f\in L^1([0, 1])\cap C((0, 1])\,$
\begin{equation}\label{HW1}
\left\{\aligned &{\phi(r) = \int_0^1 \Omega_{\frac{d-2}{2}}(rt) f(t) t^{d-1}\, dt\quad(r>0)}\quad \Rightarrow\\
 &{f(t) = \frac{1}{2^{d-2}\left[\Gamma(d/2)\right]^2}\,
 \int_0^\infty \Omega_{\frac{d-2}{2}}(rt)\, \phi(r) r^{d-1}\,dr}\quad (t>0).\endaligned\right.
\end{equation}
\end{remark}

\medskip

\subsection{Order walks}
As the radial Fourier transforms of different dimensions are known to be
interrelated by certain {\it dimension walk} transforms, the Hankel-Schoenberg transforms of
different orders turn out to be related with each other.

\smallskip

\begin{lemma}\label{lemma4} For $\,\nu>-1\,$ and $\,\alpha>0,\,\beta>0,\,$ we have
\begin{equation*}
{}_1F_2\left(\beta; \alpha +\beta,\,\nu +1; -\frac{r^2}{4}\right) =\int_0^\infty\,\Omega_\nu(rt) f(t) dt
\qquad(r\ge 0),
\end{equation*}
where $f$ is the probability density on $[0, \infty)$ defined by
\begin{equation*}
 f(t) = \frac{2}{B(\alpha, \,\beta)}\,(1-t^2)_+^{\alpha - 1}\,t^{2\beta-1}\,.
\end{equation*}
\end{lemma}

\smallskip

\begin{proof}
An elementary computation shows
\begin{align*}
\int_0^\infty t^{2k} f(t) dt = \frac{(\beta)_k}{(\alpha+\beta)_k}\,,\quad k=0, 1, 2, \cdots,
\end{align*}
and integrating termwise yields
\begin{align*}
\int_0^\infty\,\Omega_\nu(rt) f(t) dt &= \sum_{k=0}^\infty \frac{(-1)^k}{k!\,
(\nu +1)_k}\left(\frac r2\right)^{2k}\,\int_0^\infty t^{2k} f(t) dt\\
&=\sum_{k=0}^\infty \frac{(\beta)_k}{k!\,
(\alpha +\beta)_k(\nu+1)_k}\left(-\frac {r^2}{4}\right)^{k}.
\end{align*}
\end{proof}

\smallskip
Hankel-Schoenberg transforms of different orders are interrelated in the following
way, which reveals how Hankel-Schoenberg transforms generalize radial Fourier
transforms defined in \eqref{P2}.

\smallskip

\begin{theorem}\label{orderwalk}
Let $d$ be a positive integer and $\,\lambda>d/2 -1.$
\begin{itemize}
\item[\rm(i)] For each $\,r\ge 0,\,$ we have
\begin{equation*}
\Omega_{\lambda}(r) = \frac{2}{B\left(\lambda +1 -\frac d2,\,\frac d2\right)}
\int_0^\infty\Omega_{\frac{d-2}{2}}(rt)(1-t^2)_+^{\lambda -\frac d2}\, t^{d-1} dt.
\end{equation*}
\item[\rm(ii)] If $\, f\in L^1([0, \infty)),\,$ then for each $\,r\ge 0,$
\begin{align*}
&\int_0^\infty \,\Omega_{\lambda}(rt) f(t) dt
= \int_0^\infty \,\Omega_{\frac{d-2}{2}}(rt)\,I_\lambda(f)(t) t^{d-1} dt,\quad\text{where}\\
&\quad I_\lambda(f)(t) =
\frac{2}{B\left(\lambda +1 -\frac d2,\,\frac d2\right)}
\int_t^\infty\left( s^2 - t^2\right)^{\lambda - \frac d2}\,s^{-2\lambda} f(s) ds.
\end{align*}
Moreover, $\,I_\lambda(f) \in L^1\left([0, \infty), \,t^{d-1} dt\right)\,$ with
$$\int_0^\infty \left|I_\lambda(f)(t)\right| t^{d-1} dt \le \int_0^\infty \left|f(t)\right| dt.$$
\end{itemize}
\end{theorem}

\smallskip

\begin{proof}
The special choices of
$\,\nu = d/2 -1, \,\alpha = \lambda +1 - d/2, \,\beta = d/2\,$
in Lemma \ref{lemma4} gives part (i) upon noticing
\begin{equation*}
\Omega_\lambda(r) = {}_0 F_1\left(\lambda +1;\,-\frac{r^2}{4}\right).
\end{equation*}

As for part (ii), we first notice
\begin{align*}
&\int_0^\infty \int_t^\infty \left( s^2 - t^2\right)^{\lambda - \frac d2}\,s^{-2\lambda} |f(s)| ds\, t^{d-1} dt\\
&\qquad =\int_0^\infty \left[\int_0^s \left( s^2 - t^2\right)^{\lambda - \frac d2}\,t^{d-1} dt\right] s^{-2\lambda} |f(s)| ds \\
&\qquad=\int_0^\infty \left[\int_0^1 \left( 1- u^2\right)^{\lambda - \frac d2}\,u^{d-1} du\right] |f(s)| ds \\
&\qquad=\frac{B\left(\lambda +1 -\frac d2,\,\frac d2\right)}{2}\,\int_0^\infty |f(s)| ds,
\end{align*}
whence $\,I_\lambda(f) \in L^1\left([0, \infty), \,t^{d-1} dt\right)\,$ and the last estimate follows.
The stated formula is a simple consequence of part (i) on interchanging the order of integrations, which is
legitimate due to Fubini's theorem.
\end{proof}

\smallskip

\begin{remark}
If $\,d=1,\,$ part (i) reduces to Poisson's integral (P5).
The so-called descending-dimension walks of radial Fourier transforms are special
instances of this theorem. In fact,
if we take $\,\lambda = (d+k-2)/2,\,$ with $\,d, k\,$
positive integers, and write $\,I_\lambda = I_k\,$ for simplicity, then
the formula of part (ii) applied to the function $\,f(t) t^{d+k-1}\,$ yields
\begin{align}
&\int_0^\infty \,\Omega_{\frac{d+k-2}{2}}(rt) f(t) t^{d+k-1} dt
= \int_0^\infty \,\Omega_{\frac{d-2}{2}}(rt)\,I_k(f)(t) t^{d-1} dt,\nonumber\\
&\qquad I_k(f)(t) =
\frac{2}{B\left(\frac k2,\,\frac d2\right)}
\int_t^\infty\left( s^2 - t^2\right)^{\frac k2-1}\,s f(s) ds.
\end{align}
In the notation of \eqref{P2}, it reads
\begin{equation}
\F_{d+k}(f) (r) = \frac{\pi^{k/2}\,\Gamma\left(\frac d2\right)}{\Gamma\left(\frac{d+k}{2}\right)}
\,\F_{d}\left(I_k(f)\right)(r),
\end{equation}
which expresses the $(d+k)$-dimensional radial Fourier transform of $f$ as
$d$-dimensional radial Fourier transform of $I_k(f)$.
We refer to \cite{S2}, \cite{SW} and \cite{We} for more detailed results on dimension walks.
\end{remark}

\bigskip

\section{Hankel-Schoenberg transforms of binomial densities and asymptotic properties}
As the first step of our construction, we shall consider all possible binomial densities
and evaluate their Hankel-Schoenberg transforms.

\smallskip

\begin{lemma}\label{lemma6} Let $\,\lambda>-1\,$ and $\,\alpha>0,\,\beta>0.$ For each $\,r\ge 0,$
\begin{equation*}
\int_0^\infty\,\Omega_\lambda(rt) p(t) dt = {}_2F_3\left(\frac{\beta}{2},\, \frac{\beta+1}{2}; \frac{\alpha +\beta}{2},\,
\frac{\alpha +\beta+1}{2},\,\lambda +1; -\frac{r^2}{4}\right),
\end{equation*}
where $p$ is the probability density on $[0, \infty)$ defined by
\begin{equation*}
 p(t) = \frac{1}{B(\alpha, \beta)}\,(1-t)_+^{\alpha - 1}\,t^{\beta-1}\,.
\end{equation*}
\end{lemma}

\smallskip

\begin{proof}
By applying Legendre's duplication formula for the gamma function repeatedly, it is elementary to compute
\begin{align*}
\int_0^\infty t^{2k} p(t) dt = \frac{B(\alpha, 2k+\beta)}{B(\alpha, \beta)}
=\frac{\left(\frac{\beta}{2}\right)_k \left(\frac{\beta+1}{2}\right)_k}
{\left(\frac{\alpha +\beta}{2}\right)_k \left(\frac{\alpha +\beta+1}{2}\right)_k}
\end{align*}
for $\,k=0, 1, 2, \cdots,\,$ and integrating termwise yields the stated result.
\end{proof}

\smallskip

After reducing the generalized hypergeometric functions of Lemma \ref{lemma6} to the ones of type ${}_1F_2$,
we shall investigate their asymptotic properties for which
our analysis will be based on the following lemma which has been studied by many authors including
R. Askey and H. Pollard \cite{AP}, J. Steinig \cite{St} and culminated in the present form by
J. Fields and M. Ismail \cite{FI}.

\smallskip

\begin{lemma}\label{lemma7}
For $\,\rho>0,\,\nu>0,\,$ put
\begin{align*}
U(\rho, \,\nu\,; x) = {}_1 F_2\left(\nu\,;\, \rho\nu, \,\rho\nu + \frac 12\,;\,-\frac{x^2}{4}\right)
\qquad(x\in\R).
\end{align*}
\begin{itemize}
\item[\rm(i)] If $\,\rho=1,\,$ it is identical to the function $\,\Omega_{\nu-1/2}.$
\item[\rm(ii)] If $\,\rho\ne 1,\,$ then as $\,|x|\to\infty,$
\begin{align*}
&U(\rho, \,\nu\,; x) = \frac{\Gamma(2\rho\nu)}{\,\Gamma(2\rho\nu - 2\nu)\,}\,
|x|^{-2\nu}\Big[1 + O\left(|x|^{-2}\right)\Big]\\
&\quad + \frac{\Gamma(2\rho\nu)}{\,2^{\nu-1}\Gamma(\nu)\,}\,|x|^{-2\nu\left(\rho-\frac 12\right)}
\biggl[\cos\left(|x|-\rho\nu\pi + \frac{\nu\pi}{2}\right) + O\left(|x|^{-1}\right)\biggr].
\end{align*}
\item[\rm(iii)] If either $\,\rho\ge \frac 32,\,\nu>1\,$ or $\,\rho\ge 2,\,\nu>0,\,$ then
$$ U(\rho, \,\nu\,; x) \,\approx\,(1+ |x|)^{-2\nu}\quad\text{for}\quad x\in\R.$$
In particular, $\, U(\rho, \,\nu\,; x)>0\,$ for every $\,x\in\R.$
\end{itemize}
\end{lemma}

\bigskip

\section{Askey's class for $H^{\,\frac{d+1}{2}}(\R^d)$}
In the special case $\,\delta = (d+1)/2,\,$ the Bessel potential kernel $G_\delta$, which gives
a reproducing kernel for $\,H^{\,\frac{d+1}{2}}(\R^d)\,$ under the usual
inner product, coincides with the exponential of $-|\bx|$ and its Fourier transform is nothing but the Cauchy-Poisson kernel
(see appendix). To be precise, we have
\begin{equation*}
G_{\frac{d+1}{2}}(\bx) = \frac{1}{2^d \pi^{\frac{d-1}{2}}\Gamma\left(\frac{d+1}{2}\right)}
\,e^{-|\bx|}\,,\quad \widehat{G_{\frac{d+1}{2}}}(\xi) = (1+|\xi|^2)^{-\frac{d+1}{2}}.
\end{equation*}

A large class of compactly supported functions, often referred to as Askey's class (\cite{As}, \cite{G2},
\cite{We}), turn out to be also available as reproducing kernels under suitable inner
products in this case.

\smallskip

\begin{theorem}\label{askey}
For a positive integer $d$, assume that $\alpha$ satisfies $\,\alpha\ge \frac{d+1}{2}\,$ if $\,d\ge 2\,$
and $\,\alpha\ge 2\,$ if $\,d=1.\,$ Define
\begin{align*}
\Lambda_{d, \alpha}(r) = {}_1 F_2\left(\frac{d+1}{2}\,; \frac{d+\alpha+1}{2},
\frac{d+\alpha +2}{2};\, -\frac{\,r^2}{4}\right)\qquad(r\ge 0).
\end{align*}

\begin{itemize}
\item[\rm(i)] $\Lambda_{d, \alpha}$ is positive definite on $\R^d$ with
\begin{align*}
\Lambda_{d, \alpha}(r) =
\frac{1}{B(\alpha+1, \,d)}\,\int_0^\infty \Omega_{\frac{d-2}{2}}(rt)\,(1-t)_+^{\alpha}\, t^{d-1} dt\qquad(r\ge 0).
\end{align*}

\item[\rm{(ii)}] $\,0< \Lambda_{d, \alpha}(r)\le 1\,$ for each $\,r\ge 0\,$ and
\begin{align*}
\Lambda_{d, \alpha}(r) &= \frac{\Gamma(d+\alpha+1)}{\,\Gamma(\alpha)\,}\,
r^{-d-1}\Big[ 1 + O\left(r^{-2}\right)\Big]\\
 &+ \frac{\Gamma(d+\alpha+1)}{\,2^{\frac{d-1}{2}}\Gamma\left(\frac{d+1}{2}\right)} \,r^{-\frac{(d+2\alpha+1)}{2}}
\left[\cos\left(r-\frac{(d+2\alpha+1)\pi}{4}\right)+ O\left(r^{-1}\right)\right]
\end{align*}
as $\,r\to\infty.$ Moreover, $\,\Lambda_{d, \alpha}(r) \,\approx\, \left(1 + r\right)^{-d-1}\,$ for $\, r\ge 0.$

\item[\rm(iii)] $\,\Lambda_{d, \alpha}\in C([0, \infty))\cap L^1\left([0, \infty), \,r^{d-1} dr\right)\,$ and
\begin{equation*}
(1-t)_+^{\alpha} = \frac{2\Gamma(\alpha+1)\Gamma\left(\frac{d+1}{2}\right)}
{\sqrt\pi\,\Gamma(\alpha + d+1)\Gamma\left(\frac d2\right)}\,
\int_0^\infty \Omega_{\frac{d-2}{2}}(rt)\Lambda_{d, \alpha}(r) r^{d-1} dr \quad(t\ge 0).
\end{equation*}
As a consequence, the function $\,t\mapsto (1-t)^\alpha_+\,$ is positive definite on $\R^d$.
\end{itemize}
\end{theorem}

\smallskip

\begin{proof} The integral representation of part (i)
corresponds to the special case of Lemma \ref{lemma6} with $\,\lambda = \frac {d-2}{2}\,,\,\beta=d\,$
for which we replace $\alpha$ by $\alpha+1$. The positive definiteness of $\Lambda_{d, \alpha}$
is an immediate consequence of Proposition \ref{fourier}.

As to part (ii), while the uniform bound $\,\left|\Lambda_{d, \alpha}(r)\right|\le 1\,$ is a consequence of
(P1), the rest follow from Lemma \ref{lemma7} upon expressing
$$\Lambda_{d, \alpha}(r) = U\left(\frac{d+\alpha+1}{d+1}\,,\,\frac{d+1}{2}\,;\,r\right).$$

As to part (iii), the property $\,\Lambda_{d, \alpha}\in C([0, \infty))\cap L^1\left([0, \infty), \,r^{d-1} dr\right)\,$
is obvious. For $\,t>0,\,$ the stated integral representation follows by inverting the formula of part (i)
in accordance with Proposition \ref{inversion}, particularly with \eqref{HW1}. By continuity, it
continues to hold true for $\,t=0.$ Finally, the positive definiteness follows again from
Proposition \ref{fourier}.
\end{proof}

On consideration of radial extensions, we obtain the following in which
$$\gamma_{d, \alpha} = \frac{2^d \pi^{\frac{d-1}{2}}\Gamma(\alpha+1)\Gamma\left(\frac{d+1}{2}\right)}{\Gamma(\alpha+d+1)}.$$

\smallskip

\begin{corollary}\label{askey1}
For $\,\alpha\ge \frac{d+1}{2}\,$ if $\,d\ge 2\,$ and $\,\alpha\ge 2\,$ if $\,d=1,\,$ put
\begin{equation*}
A_\alpha(\bx) = (1-|\bx|)^\alpha_+\qquad(\bx\in\R^d).
\end{equation*}
Then each $A_\alpha$ is continuous and positive definite with
\begin{align*}
\widehat{A_\alpha}(\xi) = \gamma_{d, \alpha}\,\Lambda_{d, \alpha}(|\xi|)\qquad(\xi\in\R^d).
\end{align*}
As a consequence, $\,A_\alpha(\bx-\by)$ is a reproducing kernel for the Sobolev space
$\,H^{\,\frac{d+1}{2}}(\R^d)\,$ with respect to the inner product defined by
\begin{align*}
\bigl(u,\,v\bigr)_{A_\alpha(\R^d)} = (2\pi)^{-d}\cdot\frac{1}{\gamma_{d, \alpha}}\int_{\R^d}
\frac{\widehat{u}(\mathbf{\xi})\overline{\,\widehat{v}(\mathbf{\xi})} \,d\mathbf{\xi}}{
\Lambda_{d, \alpha}\left(|\mathbf{\xi}|\right)}\,.
\end{align*}
\end{corollary}

\smallskip

\begin{remark}
If we put $\,\Lambda_{d, \alpha}(\bx) = \Lambda_{d, \alpha}(|\bx|),\,\bx\in\R^d,\,$
for simplicity, then the inversion formula of part (iii) in Theorem \ref{askey} shows
$$\widehat{\Lambda_{d, \alpha}}(\xi) =
\frac{\pi^{\frac{d+1}{2}}\Gamma(\alpha +d +1)}{\Gamma(\alpha+1)\Gamma\left(\frac{d+1}{2}\right)}
\,A_\alpha(\xi)\qquad(\xi\in\R^d).$$
Thus $\Lambda_{d, \alpha}$ is an example of band-limited functions, the class of $L^2$ functions whose
Fourier transforms are compactly supported.
\end{remark}

In the odd dimensional case, $\Lambda_{d, \alpha}$ is expressible in terms of algebraic and trigonometric functions
if $\alpha$ happens to be an integer. As illustrations, we present the following examples:

\begin{itemize}
\item[(a)] In the case $\,d=1,$ the formula of part (i) in Theorem \ref{askey} reduces to
\begin{equation*}
\Lambda_{1, \alpha}(r) = (\alpha+1)\,\int_0^1 \cos(xt)(1-t)^\alpha dt
\quad(\alpha\ge 2).
\end{equation*}

With the choice of minimal $\,\alpha =2\,$ and $\,\alpha=3,\,$ we have
\begin{align*}
\Lambda_{1, 2}(r) &= \frac{6}{\,r^2\,}\,\left( 1- \frac{\sin r}{r}\right)\,,\\
\Lambda_{1, 3}(r) &= \frac{12}{\,r^2\,}\left\{ 1- \left[\frac{\sin (r/2)}{r/2}\right]^2\right\}
\end{align*}
in which each formula must be understood as the limiting value at $\,r=0\,$.

\item[(b)] In the case $\,d=3,$ the formula of part (i) in Theorem \ref{askey} reduces to
\begin{equation*}
\Lambda_{3, \alpha}(r) = \frac{(\alpha+3)(\alpha+2)(\alpha+1)}{2r}\,\int_0^1 \sin(rt)(1-t)^\alpha t dt
\quad(\alpha\ge 2).
\end{equation*}

With the choice of minimal $\,\alpha = 2,\,$ we have
\begin{align*}
\Lambda_{3, 2}(r) = \frac{60}{\,r^4\,} \left( 2 + \cos r -
3 \,\frac{\sin r}{r}\right)
\end{align*}
with the same interpretation at $\,r=0\,$ as above.
\end{itemize}

\bigskip

\section{Compactly supported reproducing kernels for $H^{\,\delta}(\R^d)$
with $\,\delta>\max\, (1, \, d/2)$}
Due to an obvious cancellation effect, the generalized hypergeometric function of Lemma \ref{lemma6}
in the special case $\,\beta = 2\lambda +1\,$ reduces to
\begin{align}\label{W0}
&{}_1F_2\left(\frac{2\lambda +1}{2}; \frac{\alpha + 2\lambda +1}{2},\,
\frac{\alpha + 2\lambda +2}{2}; -\frac{r^2}{4}\right)\nonumber\\
&\qquad\qquad = \frac{1}{B(\alpha, 2\lambda +1)}\int_0^\infty\,\Omega_\lambda(rt)(1-t)_+^{\alpha - 1}\,t^{2\lambda}dt.
\end{align}

Expressing in the form of $U$-function defined in Lemma \ref{lemma7}, it is
simple to find that this function is strictly positive if $\,\lambda>1/2,\,\alpha\ge \lambda + 1/2.\,$
The choice of minimal value $\,\alpha = \lambda + 1/2\,$ leads to
\begin{align}\label{W1}
&{}_1F_2\left(\lambda + \frac 12 ; \frac 32\left(\lambda + \frac 12\right),\, \frac 32\left(\lambda + \frac 12\right) + \frac 12;
 -\frac{r^2}{4}\right)\nonumber\\
&\qquad\qquad = \frac{1}{B\left(\lambda + \frac 12, \,2\lambda +1\right)}\int_0^\infty\,\Omega_\lambda(rt)(1-t)_+^{\lambda - \frac 12}\,t^{2\lambda}dt.
\end{align}

Rearranging parameters $\,\lambda + 1/2 = \delta\,$ and representing the last Hankel-Schoenberg transforms
in terms of radial Fourier transforms, that is, those integrals with kernels $\Omega_{\frac{d-2}{2}}$,
we are led to the following class of functions.

\smallskip

\begin{definition}
For a positive integer $d$ and $\,\delta>d/2,\,$ define
\begin{align*}
\Phi_{d, \delta}(t) = \frac{1}{B\left(2\delta -d,\,\delta\right)}
\int_t^1 (s^2- t^2)^{\delta -\frac {d+1}{2}}\,
(1-s)^{\delta -1}\,ds
\end{align*}
for $\,0\le t\le 1\,$ and zero otherwise.
\end{definition}

\smallskip

\begin{lemma}\label{lemma8}
For a positive integer $d$ and $\,\delta>d/2,\,$ the integral in the definition of $\Phi_{d, \delta}$
converges and the following properties hold:
\begin{itemize}
\item[\rm(i)] $\Phi_{d, \delta}$ is continuous, strictly decreasing on $[0, 1]$ and
$\,0\le \Phi_{d, \delta}\le 1.$

\item[\rm(ii)] $\,\Phi_{d, \delta}(t)\,\approx\, (1-t)^{2\delta  - \frac{d+1}{2}}\,$ on $[0, 1].$

\item[\rm(iii)] If $\,\delta = \frac{d+1}{2}\,,\,$ then
$\,\Phi_{d, \delta} (t) = (1- t)_+^{\frac{d+1}{2}}\,.$

\item[\rm(iv)] If $\,\delta>\frac{d+1}{2}\,,\,$ then for $\,0\le t\le 1,$
\begin{align*}
\Phi_{d, \delta}(t) = \frac{1}{B\left(2\delta -d-1,\,\delta+1\right)}
\int_t^1 (s^2- t^2)^{\delta -\frac {d+3}{2}}\,
s(1-s)^{\delta}\,ds.
\end{align*}
\end{itemize}
\end{lemma}

\smallskip

\begin{proof}
For $\,\delta\ge \frac{d+1}{2}\,,\,$ as the function $\Phi_{d, \delta}$ is dominated by
$$
\frac{1}{B\left(2\delta -d,\,\delta\right)}
\int_0^1 s^{2\delta -d -1} (1-s)^{\delta -1}\,ds = 1,
$$
the convergence of the defining integral is obvious.
Under the transformation $\,s\mapsto \theta + (1-\theta) t,\,\,0 \le\theta\le 1, \,$ we may write
\begin{align*}
&\qquad\Phi_{d, \delta}(t) = (1-t)^{2\delta -\frac{d+1}{2}}\,V(t),\quad\text{where}\\
& V(t) = \frac{1}{B\left(2\delta -d, \,\delta\right)}
\int_0^1 \theta^{\delta - \frac{d+1}{2}}(1-\theta)^{\delta-1}\big[2t + \theta(1-t)\big]^{\delta- \frac{d+1}{2}}
d\theta\,.
\end{align*}
In the case $\, \frac{d}{2}<\delta<\frac{d+1}{2}\,,$ if we observe
$$ 2^{\delta- \frac{d+1}{2}}\le \big[2t + \theta(1-t)\big]^{\delta- \frac{d+1}{2}}
\le \theta^{\delta- \frac{d+1}{2}}$$
for $\,0\le t\le 1\,$ and for each fixed $\,\theta>0,$ it is simple to infer that
$V(t)$ converges uniformly on $[0, 1]$ with
$\,0\le V(t)\le 1\,$ and hence $\Phi_{d, \delta}$ is well defined. Bounding $V(t)$ in this way,
we also deduce part (ii) plainly.

As the convergence is ensured, part (i) can be verified easily. Part (iii) is trivial and
part (iv) is a simple consequence of integrating by parts.
\end{proof}

\smallskip

\begin{remark} Noteworthy are the following special instances of part (iv).

\begin{itemize}
\item [(a)] In the case $\,\delta = d/2 + k + 1/2,\,k\in\mathbb{N},\,$
$\Phi_{d, \delta}$ coincides with Wendland's function $P_{d, k}$, defined in \eqref{G3}, in the odd dimensions.
\item[(b)] In the case $\,\delta = d/2 + m + 1,\,$ with $m$ a nonnegative integer,
$\Phi_{d, \delta}$ coincides with Schaback's function $R_{d, m}$, defined in \eqref{G4}, in every dimension.
Likewise, if $\,\delta = (d+1)/2 + \alpha,\,\alpha>0,\,$ $\Phi_{d, \delta}$ coincides with
the function $S_{d, \alpha}$ of Chernih and Hubbert, defined in \eqref{G5}, in every dimension.
\end{itemize}
\end{remark}

\smallskip

In the statement below, we shall denote
\begin{equation}
\omega_{d, \delta} =\frac{ 2^{1 -d}\,\Gamma\left(3\delta\right)\Gamma\left(\delta -\frac d2\right)}
{\Gamma\left(\delta\right)\Gamma\left(3\delta-d\right)\Gamma\left(\frac d2\right)}\,.
\end{equation}

\smallskip

\begin{theorem}\label{wend}
For a positive integer $d$ and $\,\delta>\max\left(1, \,d/2\right),$ define
\begin{align*}
W_\delta(r) = {}_1 F_2\left(\delta\,; \frac {3\delta}{2}, \,\frac {3\delta +1}{2}\,;\, -\frac{\,r^2}{4}\right)
\qquad(r\ge 0).
\end{align*}

\begin{itemize}
\item[\rm(i)] $W_\delta$ is positive definite on $\R^d$ with
\begin{align*}
W_\delta(r) &=
\frac{1}{B\left(\delta, \,2\delta\right)}\int_0^\infty \Omega_{\delta- \frac 12}(rt)
(1-t)_+^{\delta-1}\, t^{2\delta -1}\,dt\\
&= \omega_{d, \delta} \int_0^\infty \Omega_{\frac{d-2}{2}}(rt)\,\Phi_{d, \delta}(t)\, t^{d-1} dt.
\end{align*}

\item[\rm{(ii)}] $\,0<W_\delta(r)\le 1\,$ for each $\,r\ge 0\,$ and as $\,r\to\infty,$
\begin{align*}
W_\delta(r) &= \frac{\Gamma\left(3\delta\right)}{\Gamma\left(\delta\right)}
\,r^{-2\delta}\biggl[ 1 + \frac{\cos\left(r -\delta\pi\right)}{2^{\delta -1}}\biggr]
+ O\left(r^{-2\delta -1}\right)\,.
\end{align*}
Moreover, $\,W_\delta(r) \,\approx\, \left(1 + r\right)^{-2\delta}\,$ for $\,r\ge 0.\,$

\item[\rm(iii)] $\,W_\delta\in C([0, \infty))\cap L^1\left([0, \infty), \,r^{d-1} dr\right)\,$ and
\begin{equation*}
\Phi_{d, \delta} (t) = \frac{1}{2^{d-2}\,\left[\Gamma(d/2)\right]^2\,\omega_{d, \delta}}\,
\int_0^\infty \Omega_{\frac{d-2}{2}}(rt) W_\delta(r) r^{d-1} dr \qquad(t\ge 0).
\end{equation*}
As a consequence, $\Phi_{d, \delta}$ is positive definite on $\R^d$.
\end{itemize}
\end{theorem}

\smallskip

\begin{proof}
If we set $\,\lambda = \delta - 1/2\,$ in the representation \eqref{W1}, we obtain
\begin{align*}
W_\delta(r) &=
\frac{1}{B\left(\delta, \,2\delta\right)}\int_0^\infty \Omega_{\delta- \frac 12}(rt)
(1-t)_+^{\delta-1}\, t^{2\delta -1}\,dt\\
&= C(d, \delta)\int_0^\infty \Omega_{\frac{d-2}{2}}(rt)\,\Phi_{d, \delta}(t)\, t^{d-1} dt,
\end{align*}
where the latter follows by the order-walk transform of Theorem \ref{orderwalk} and
$$C(d, \delta) = \frac{B\left(2\delta - d,\,\delta\right)}
{B\left(\delta +\frac 12 - \frac d2, \,\frac d2\right) B(\delta, \,2\delta)}\,.$$
Simplifying with the aid of Legendre's duplication formula for the gamma function, it is elementary to
see $\,C(d, \delta) = \omega_{d, \delta}.$ The positive definiteness of $W_\delta$
is an immediate consequence of Proposition \ref{fourier} and part (i) is proved.

In view of the identification
$$W_\delta(r) = U\left(\frac 32,\,\delta\,;\,r\right),$$
part (ii) is a consequence of Lemma \ref{lemma7} and Lemma \ref{lemma8}.

As to part (iii), that $\,W_\delta\in C([0, \infty))\cap L^1\left([0, \infty), \,r^{d-1} dr\right)\,$
is obvious. For $\,t>0,\,$ the stated representation follows by inverting the formula of part (i)
in accordance with \eqref{HW1}. By continuity, it
continues to hold true for $\,t=0.$ Finally, the positive definiteness follows again from
Proposition \ref{fourier}.
\end{proof}

As an immediate corollary, we obtain what we aim to accomplish. To simplify notation, we shall
write
\begin{equation}
\zeta_{\,d, \delta} =\frac{2\pi^{d/2}}{\Gamma(d/2)}\cdot \frac{1}{\omega_{d, \delta}}.
\end{equation}

\smallskip

\begin{corollary}\label{wend1}
For $\,\delta>\max\left(1, \,d/2\right),$ let $\,\Phi_{d, \delta} (\bx) = \Phi_{d, \delta} (|\bx|),\,\bx\in\R^d.$
Then $\Phi_{d, \delta}$ is continuous and positive definite with
\begin{align*}
\widehat{\Phi_{d, \delta}}(\xi) = \zeta_{\,d, \delta}\,W_\delta(|\xi|)\qquad(\xi\in\R^d).
\end{align*}
As a consequence, $\,\Phi_{d, \delta}(\bx-\by)$ is a reproducing kernel for the Sobolev space
$\,H^{\,\delta}(\R^d)\,$ with respect to the inner product defined by
\begin{align*}
\bigl(u,\,v\bigr)_{\Phi_{d, \delta}(\R^d)} = (2\pi)^{-d}\cdot\frac{1}{\zeta_{\,d, \delta}}\,\int_{\R^d}
\frac{\widehat{u}(\mathbf{\xi})\overline{\,\widehat{v}(\mathbf{\xi})} \,d\mathbf{\xi}}{
W_\delta\left(|\mathbf{\xi}|\right)}\,.
\end{align*}
\end{corollary}

\begin{remark} In the special case $\,\delta = (d+1)/2 +\alpha, \,\alpha>0,\,$  A. Chernih and S. Hubbert
also obtained the Fourier transform $W_\delta$ (Theorem 2.1, \cite{CH}), but the authors did not give the integral representation
formula nor the inversion formula as stated in the first equation of part (ii), part (iii) of Theorem \ref{wend}, respectively.
We supplement a few computational aspects as follows.

\begin{itemize}
\item[(a)] In some special instances, it is possible to evaluate $\Phi_{d, \delta}$
in closed forms as the following list shows.

\bigskip
\bigskip
\begin{tabular}{|l|lc|} \hline
 & \multicolumn{2}{|c|}{\textbf{$\Phi_{d, \delta}(r)$ on the interval $\,[0, 1]\,$}}\\\hline
$\,\,\delta = \frac {d+1}{2}\,\,$ & $d\ge 2$   & $ (1-r)^{\frac{d+1}{2}}$  \\ \hline
& $d=1$ &   $\,(1-r)^3 \,(1+3r)\,$   \\ \cline{2-3}
$\,\,\delta = 2\,\,$ &$d=2$ &  $\,\left( 1+ 2r^2\right)\sqrt{1-r^2} -3r^2
\log \left(\frac{ 1+\sqrt{1-r^2}}{r}\right)\,$\\ \cline{2-3}
& $d=3$ &$ (1- r)^2\,$  \\ \hline
& $ d=1$ &$\,(1-r)^5 ( 1-2r + 8r^2)$\\ \cline{2-3}
$\,\,\delta = 3\,\,$ & $ d=2$ & $ \frac 14 \Big[(4- 28 r^2 - 81 r^4)\sqrt{1-r^2}
$\\
& & $ \qquad\qquad \qquad\quad+\,\,15 r^4 (6 + r^2)\log \left(\frac{ 1+\sqrt{1-r^2}}{r}\right)\Big]\,$    \\ \cline{2-3}
&$ d=3$ & $ (1-r)^4 (1+4r) $\\ \hline
\end{tabular}
\bigskip
\bigskip

\item[(b)] In the case when $\delta$ is an integer, one may use the representation formula of
part (i), Theorem \ref{wend}, to express $W_\delta$ in a closed form
involving algebraic and trigonometric functions. To illustrate,
let us take
$$W_{2}(r) = {}_1F_2\left(2; 3, \,\frac 72; -\frac{r^2}{4}\right)\qquad(r\ge 0).$$

We evaluate
\begin{align*}
W_{2}(r) &= 20\int_0^\infty\Omega_{3/2}(rt)(1-t)_+\,t^3\,dt\\
&=\frac{60}{r^3}\left\{\int_0^1\sin(rt) (1-t) dt - r\int_0^1\cos(rt) (1-t) tdt\right\}\\
&= \frac{120}{r^4}\,\left( 1+ \frac{\cos r}{2}\right) - \frac{180\, \sin r}{r^5}\,.
\end{align*}
We should point out this closed form is consistent with the asymptotic formula
stated in part (ii) of Theorem \ref{wend} which reads
$$W_{2}(r) =\frac{120}{r^4}\,\left( 1+ \frac{\cos r}{2}\right) + O\left(r^{-5}\right).$$
\end{itemize}
\end{remark}

\bigskip

\section{A smoother family of compactly supported reproducing kernels}
Due to the restriction $\,\delta>\max\,(1, \,d/2),\,$ there are missing cases in the preceding
results, namely, the cases $\,1/2<\delta\le 1\,$ for the one-dimensional Sobolev spaces $H^\delta(\R).$
Although the particular instance $\,\delta =1\,$ is covered in Corollary \ref{askey1}, the case
$\,1/2<\delta<1\,$ is still left out.

The purpose of this section is to provide compactly supported reproducing kernels
in such missing cases. As a matter of fact, we shall construct another class of compactly supported reproducing kernels
which suit to the Sobolev spaces $H^\delta(\R^d)$ of any order $\,\delta>d/2\,$ without any restriction.

The key idea is to exploit the lemma of J. Fields and M. Ismail, Lemma \ref{lemma7}, in such a way that
the strict positivity of the generalized hypergeometric function of \eqref{W0} is assured in the range $\,\lambda + 1/2>0,\,\alpha\ge 2\lambda+1.\,$
Choosing the minimal $\,\alpha= 2\lambda+1\,$ and setting  $\,\lambda + 1/2 = \delta,\,$
it reduces to
\begin{align}\label{S0}
&{}_1F_2\left(\delta; 2\delta, \,2\delta + \frac 12; -\frac{r^2}{4}\right)\nonumber\\
&\qquad = \frac{1}{B(2\delta, \,2\delta)}\int_0^\infty\,\Omega_{\delta- \frac 12}(rt)(1-t)_+^{2\delta -1}\,t^{2\delta-1}dt,
\end{align}
which is strictly positive for any $\,\delta>0.$

For $\,\delta>d/2,\,$ an application of order-walk transformation yields
\begin{align}\label{S1}
\int_0^\infty\,\Omega_{\delta- \frac 12}(rt)(1-t)_+^{2\delta -1}\,t^{2\delta-1}dt
= \int_0^\infty\,\Omega_{\frac{d-2}{2}}(rt) I_\delta(t) t^{d-1} dt
\end{align}
in which $I_\delta$ stands for the function supported in $[0, 1]$ and defined by
\begin{align*}
I_\delta(t) = \frac{2}{B\left(\delta + \frac 12 - \frac d2,\,\frac d2\right)}\,\int_t^1 (s^2 - t^2)^{\delta - \frac{d+1}{2}}
(1-s)^{2\delta-1}\,ds
\end{align*}
for $\,0\le t\le 1.$ Normalizing the constant, we introduce

\smallskip

\begin{definition}
For a positive integer $d$ and $\,\delta>d/2,\,$ define
\begin{align*}
\Psi_{d, \delta}(t) = \frac{1}{B\left(2\delta -d,\,2\delta\right)}
\int_t^1 (s^2- t^2)^{\delta -\frac {d+1}{2}}\,
(1-s)^{2\delta -1}\,ds
\end{align*}
for $\,0\le t\le 1\,$ and zero otherwise.
\end{definition}

\smallskip

Being of similar nature with $\Phi_{d, \delta}$, we deduce its basic properties in the same way
as stated and proved in Lemma \ref{lemma8}.

\smallskip

\begin{lemma}\label{lemma9}
For a positive integer $d$ and $\,\delta>d/2,\,$ the integral in the definition of $\Psi_{d, \delta}$
converges and the following properties hold:
\begin{itemize}
\item[\rm(i)] $\Psi_{d, \delta}$ is continuous, strictly decreasing on $[0, 1]$ and
$\,0\le \Psi_{d, \delta}\le 1.$

\item[\rm(ii)] $\,\Psi_{d, \delta}(t)\,\approx\, (1-t)^{3\delta  - \frac{d+1}{2}}\,$ on $[0, 1].$

\item[\rm(iii)] If $\,\delta = \frac{d+1}{2}\,,\,$ then
$\,\Psi_{d, \delta} (t) = (1- t)_+^{d+1}\,.$

\item[\rm(iv)] If $\,\delta>\frac{d+1}{2}\,,\,$ then for $\,0\le t\le 1,$
\begin{align*}
\Psi_{d, \delta}(t) = \frac{1}{B\left(2\delta -d-1,\,2\delta+1\right)}
\int_t^1 (s^2- t^2)^{\delta -\frac {d+3}{2}}\,
s(1-s)^{2\delta}\,ds.
\end{align*}
\end{itemize}
\end{lemma}

\smallskip

Combining \eqref{S0}, \eqref{S1} in terms of $\Psi_{d, \delta}$, we obtain the following analog
of Theorem \ref{wend} without difficulty in which we write
\begin{equation}
\tau_{d, \delta} =\frac{ 2^{1-d}\,\Gamma\left(4\delta\right)\Gamma\left(\delta -\frac d2\right)}
{\Gamma\left(\delta\right)\Gamma\left(4\delta-d\right)\Gamma\left(\frac d2\right)}\,.
\end{equation}

\smallskip

\begin{theorem}\label{smooth}
For a positive integer $d$ and $\,\delta> d/2,\,$ define
\begin{align*}
Q_\delta(r) = {}_1 F_2\left(\delta\,;  2\delta, \,2\delta + \frac 12\,;\, -\frac{\,r^2}{4}\right)
\qquad(r\ge 0).
\end{align*}

\begin{itemize}
\item[\rm(i)] $Q_\delta$ is positive definite on $\R^d$ with
\begin{align*}
Q_\delta(r) &=
\frac{1}{B\left(2\delta, \,2\delta\right)}\int_0^\infty \Omega_{\delta- \frac 12}(rt)
(1-t)_+^{2\delta-1}\, t^{2\delta -1}\,dt\\
&= \tau_{d, \delta} \int_0^\infty \Omega_{\frac{d-2}{2}}(rt)\,\Psi_{d, \delta}(t)\, t^{d-1} dt.
\end{align*}

\item[\rm{(ii)}] $\,0<Q_\delta(r)\le 1\,$ for each $\,r\ge 0\,$ and as $\,r\to\infty,$
\begin{align*}
Q_\delta(r) &= \frac{\Gamma\left(4\delta\right)}{\Gamma\left(2\delta\right)}
\,r^{-2\delta} + O\left(r^{-\,\min\,(2\delta +2,\,3\delta)}\right)\,.
\end{align*}
Moreover, $\,Q_\delta(r) \,\approx\, \left(1 + r\right)^{-2\delta}\,$ for $\,r\ge 0.\,$

\item[\rm(iii)] $\,Q_\delta\in C([0, \infty))\cap L^1\left([0, \infty), \,r^{d-1} dr\right)\,$ and
\begin{equation*}
\Psi_{d, \delta} (t) = \frac{1}{2^{d-2}\,\left[\Gamma(d/2)\right]^2\,\tau_{d, \delta}}\,
\int_0^\infty \Omega_{\frac{d-2}{2}}(rt) Q_\delta(r) r^{d-1} dr \qquad(t\ge 0).
\end{equation*}
As a consequence, $\Psi_{d, \delta}$ is positive definite on $\R^d$.
\end{itemize}
\end{theorem}

\smallskip

As an immediate corollary, we obtain the following in which
\begin{equation}
\eta_{\,d, \delta} =\frac{2\pi^{d/2}}{\Gamma(d/2)}\cdot \frac{1}{\tau_{d, \delta}}.
\end{equation}

\smallskip

\begin{corollary}\label{smooth1}
For $\,\delta> d/2,$ let $\,\Psi_{d, \delta} (\bx) = \Psi_{d, \delta} (|\bx|),\,\bx\in\R^d.$
Then $\Psi_{d, \delta}$ is continuous and positive definite with
\begin{align*}
\widehat{\Psi_{d, \delta}}(\xi) = \eta_{\,d, \delta}\,Q_\delta(|\xi|)\qquad(\xi\in\R^d).
\end{align*}
As a consequence, $\,\Psi_{d, \delta}(\bx-\by)$ is a reproducing kernel for the Sobolev space
$\,H^{\,\delta}(\R^d)\,$ with respect to the inner product defined by
\begin{align*}
\bigl(u,\,v\bigr)_{\Psi_{d, \delta}(\R^d)} = (2\pi)^{-d}\cdot\frac{1}{\eta_{\,d, \delta}}\,\int_{\R^d}
\frac{\widehat{u}(\mathbf{\xi})\overline{\,\widehat{v}(\mathbf{\xi})} \,d\mathbf{\xi}}{
Q_\delta\left(|\mathbf{\xi}|\right)}\,.
\end{align*}
\end{corollary}

\smallskip

\begin{remark} In view of parts (ii), (iii) of Lemma \ref{lemma8}, it is evident that
$\Psi_{d, \delta}$ is much smoother than $\Phi_{d, \delta}$ is, if both parameters $\,d, \delta\,$ are
fixed. A possible disadvantage in practical applications, however, is that $\Psi_{d, \delta}$
involves higher algebraic powers than $\Phi_{d, \delta}$ does.

\begin{itemize}
\item[(a)] As illustrations, we have the following evaluations:

\bigskip
\bigskip
\begin{tabular}{|l|lc|} \hline
 & \multicolumn{2}{|c|}{\textbf{$\Psi_{d, \delta}(r)$ on the interval $\,[0, 1]\,$}}\\\hline
$d\ge 1$  & $\delta = \frac {d+1}{2}$  & $ (1-r)^{d+1}$  \\ \hline

& $\delta = \frac 32$ &   $\,\frac 12 \Big[( 2 + 13 r^2) \sqrt{ 1- r^2} -3 r^2 ( 4 + r^2) \log \left( \frac{1+ \sqrt{1- r^2}}{r}\right)
\Big]$   \\ \cline{2-3}
$d=1$ &$\delta = 2$ &  $(1- r)^5 (1+5r)$\\ \cline{2-3}
& $\delta =3$ &$ (1- r)^8 ( 1+ 8r + 21 r^2)$  \\ \hline

& $\delta = \frac 52$ &   $(1-r)^6 (1+6r)\,$   \\ \cline{2-3}
$d=2$ &$\delta = \frac 72$ &  $\frac 13 (1- r)^9 (3 + 27 r + 80 r^2)$\\ \cline{2-3}
& $\delta = \frac 92$ &$ (1- r)^{12} (1 + 12 r + 57 r^2 + 112 r^3)$  \\ \hline

$d=3$ & $ \delta =3$ & $ (1- r)^7 ( 1+ 7r)$    \\ \cline{2-3}
&$ \delta = 4$ & $ (1-r)^{10} ( 1 + 10 r + 33 r^2) $\\ \hline
\end{tabular}
\bigskip
\bigskip


\item[(b)] As before, one may use the representation formula of
part (i), Theorem \ref{smooth}, to express $Q_\delta$ in a closed form
involving algebraic and trigonometric functions in some instances. To illustrate,
let us take
$$Q_{2}(r) = {}_1F_2\left(2; 4, \,\frac 92; -\frac{r^2}{4}\right)\qquad(r\ge 0).$$

In this special case, we evaluate
\begin{align*}
Q_{2}(r) &= 140\int_0^\infty\Omega_{3/2}(rt)(1-t)_+^3\,t^3\,dt\\
&=\frac{420}{r^3}\left\{\int_0^1\sin(rt) (1-t)^3 dt - r\int_0^1\cos(rt) t(1-t)^3 dt\right\}\\
&= \frac{840}{r^7}\,\left( r^3  - 12 r + 15 \sin r - 3r \cos r\right),
\end{align*}
which is consistent with the asymptotic formula
$$Q_{2}(r) =\frac{840}{r^4} \,+ O\left(r^{-6}\right).$$
\end{itemize}
\end{remark}

\section{Appendix: Bessel potential kernels}
In addition to the Fourier transform formulas, the Bessel potential kernels $G_\delta$
(or Mat\'ern functions) possess a number of important properties and arise
in many areas of Mathematics with various disguises. As we are concerned with constructing
possible replacements of Bessel potential kernels in the subject of reproducing kernels for Sobolev spaces,
it may be instructive to recall some of their very basic properties (see \cite{AS}, \cite{Wa}).

\begin{itemize}
\item [(a)] Each $G_\delta$ is smooth away from the origin and subject to the asymptotic
behavior, modulo multiplicative constants, described as follows:
\begin{align*}
&{\rm (i)}\quad\text{As}\quad |\bx|\to \infty,\quad
G_\delta(\bx)\, \sim\, e^{-|\bx|}\,|\bx|^{\delta - \frac{d+1}{2}}\,.\\
&{\rm (ii)}\quad\text{As}\quad |\bx|\to 0,\quad\,\,
G_\delta(\bx)\, \sim\,\left\{\begin{aligned}
&{\quad\,\, 1} &{\quad\text{if}\quad \delta> d/2\,},\\
&{-\log |\bx|} &{\quad\text{if}\quad \delta = d/2\,},\\
&{\quad |\bx|^{2\delta-d}} &{\quad\text{if}\quad \delta< d/2\,}.\end{aligned}\right.
\end{align*}

\item[(b)] Due to Schl\"afli's integral representations,
\begin{align}\label{K3}
K_{\alpha}(z) &= \frac{\sqrt{\pi}}{\Gamma(\alpha + 1/2)}\,\left(\frac{z}{2}\right)^{\alpha}
\int_{1}^{\infty} e^{-zt}\left(t^{2}-1\right)^{\alpha-\frac{1}{2}}\,dt\nonumber\\
&=\sqrt{\frac{\pi}{2}\,} \frac{e^{-z} z^\alpha}{\Gamma(\alpha+1/2)}\,
\int_{0}^{\infty} e^{- zt} \left[ t\left( 1 + \frac t2\right)\right]^{\alpha - \frac 12}\,dt,
\end{align}
which is valid for $\,\alpha>-1/2\,$ and $\,z>0,$ it is easy to see
$$ G_{\frac{d+1}{2}}(\bx) = \frac{K_{\frac 12}(|\bx|)\sqrt{|\bx|}}{2^{d-\frac 12}\,\pi^{\frac d2}\Gamma\left(\frac{d+1}{2}\right)} =
\frac{e^{-|\bx|}}{2^{d}\,\pi^{\frac {d-1}{2}}\Gamma\left(\frac{d+1}{2}\right)}\,.$$

\item[(c)] More generally, if $m$ is a nonnegative integer, then
\begin{equation}
G_{m+ \frac{d+1}{2}}(\bx) = \frac{e^{-|\bx|}\,|\bx|^m}{2^{m+d}\,\pi^{\frac{d-1}{2}}\,\Gamma\left(m+ \frac{d+1}{2}\right)}
\sum_{k=0}^m\frac{(m+k)!\,}{k! (m-k)!} \,(2|\bx|)^{-k}\,,
\end{equation}
which can be deduced easily from Schl\"afli's integrals.

\end{itemize}

\end{document}